\documentclass{amsart}

 \usepackage{graphicx,color}
\usepackage{geometry}
\usepackage{color, graphicx}
\usepackage[all]{xy}

\usepackage{calrsfs}
\usepackage{times}
\usepackage[scaled=0.92]{helvet}

\theoremstyle{plain}
\newtheorem{introprop}{Proposition}
\newtheorem{introthm}[introprop]{Theorem}

\theoremstyle{definition}
\newtheorem*{introrem}{Remark}
\swapnumbers
\theoremstyle{plain}
\newtheorem{thm}{Theorem}[section]
\newtheorem{prop}[thm]{Proposition}

\newtheorem{lem}[thm]{Lemma}
\newtheorem{lemma}[thm]{Lemma}
\theoremstyle{definition}
\newtheorem{df}[thm]{Definition}

\newtheorem{rem}[thm]{Remark}
\newtheorem{remark}[thm]{Remark}

\newtheorem{se}[thm]{}

\usepackage{amsfonts}
\usepackage{amssymb}

\renewcommand{\tilde}{\widetilde}
\renewcommand{\bar}{\overline}
\DeclareMathOperator{\fhi}{\varphi}
\DeclareMathOperator{\Z}{\mathbf{Z}}
\DeclareMathOperator{\Hom}{\mathrm{Hom}}
\DeclareMathOperator{\Aut}{\mathrm{Aut}}

\DeclareMathOperator{\kar}{\mathrm{char}}
\DeclareMathOperator{\tr}{\mathrm{tr}}

\DeclareMathOperator{\id}{\mathrm{id}}
\DeclareMathOperator{\maks}{\mathfrak{m}}

\DeclareMathOperator{\kate}{\mathsf{Art}_{\mathit k}}%{\mathcal{C}_{\mathit k}}
\DeclareMathOperator{\katez}{\widehat{\mathsf{Art}}_{\mathit k}}%\mathcal{C}}_{\mathit k}}

\DeclareMathOperator{\psl}{\![\![\!}
\DeclareMathOperator{\psr}{\!]\!]}

\DeclareMathOperator{\PB}{\mathfrak{B}}
\DeclareMathOperator{\PGL}{\mathrm{PGL}}
\DeclareMathOperator{\GL}{\mathrm{GL}}
\renewcommand{\geq}{\geqslant}
\renewcommand{\leq}{\leqslant}

\newlength{\myVSpace}% the height of the box
\newcommand\xstrut{\raisebox{-.5\myVSpace}% symmetric behaviour, 
  {\rule{0pt}{\myVSpace}}%
}
\begin{document}

\date{\today\ (version 1.2)} 
\title[universal deformations]{Which weakly ramified group actions admit\\  a universal formal deformation?}
\author[J.~Byszewski]{Jakub Byszewski}
\author[G.~Cornelissen]{Gunther Cornelissen}
\address{Mathematisch Instituut, Universiteit Utrecht, Postbus 80.010, 3508 TA Utrecht, Nederland}
\email{\{j.j.byszewski,g.cornelissen\}@uu.nl}

\thanks{Research done in part at the Universit\'e de Paris-Sud XI, supported by the European Community under the contract A.A.G.\ (Arithmetic Algebraic Geometry). J.B.\ thanks the Universit\'e de Paris-Sud XI and the Universit\'e de Versailles Saint-Quentin for their hospitality, in particular \'Etienne Fouvry and Ariane M\'ezard. We thank Ariane M\'ezard for many useful comments on the text.} 
%\address{Mathematisch Instituut, Universiteit Utrecht, Postbus 80.010, 3508 TA Utrecht, Nederland}
%

\begin{abstract} \noindent Consider a formal (mixed-characteristic) deformation functor $D$ of a representation of a finite group $G$ as automorphisms of a power series ring $k\psl t \psr$ over a perfect field $k$ of positive characteristic. Assume that the action of $G$ is weakly ramified, i.e., the second ramification group is trivial. Examples of such representations are provided by a group action on an ordinary curve: the action of a ramification group on the completed local ring of any point on such a curve is weakly ramified. 

We prove that the only such $D$ that are not pro-representable occur if $k$ has characteristic two and $G$ is of order two or isomorphic to a Klein group. Furthermore, we show that only the first of those has a non-pro-representable equicharacteristic deformation functor. 

\end{abstract}

\maketitle

\section*{Introduction}
The optimal situation in deformation theory occurs when a universal object exists --- when a deformation functor is (pro-)representable. For example, this happens in the formal deformation theory of a group action on a projective curve of genus $g \geq 2$ (\cite{BM}), or for absolutely irreducible Galois representations (\cite{Mazur}) --- The latter example played a decisive r\^{o}le in the proof of Fermat's Last Theorem. Equally often, one doesn't expect or cannot establish (pro-)representability, and the remedy is the construction of a so called ``versal hull'' for the deformation functor (\cite{Schlessinger}). This is the classical approach to the local version of the first example: the  action of a finite group on the completed local ring of a point of a curve. In this work, we will prove that some of these versal hulls are actually universal, though not by a standard method (unfortunately, the literature seems to be marred by the use of the expression ``is not (pro-)representable'' instead of ``has not been established to be (pro-)representable by this or that method'', which leads to a lot of confusion). The problem came up naturally in dealing with d\'evissage (\cite{Kuba}) and in clarifying some points in the computation of versal (!) hulls in \cite{CM} (cf. Remark \ref{zadziornosc} infra: the universality for $G=\Z/p$ is used to compute the versal ring for general $G$).

First, we set up the precise notation to explain the results.  Denote by $k$ a fixed perfect field of positive characteristic $p$, and by $W(k)$ the ring of Witt vectors of $k$. Let $\kate$ denote the category of local artinian $W(k)$-algebras with residue field $k$ and local morphisms of $W(k)$-algebras; and let $\katez$ denote the category of complete local noetherian $W(k)$-algebras with residue field $k$ and local morphisms of $W(k)$-algebras.
Then $\kate$ is a full subcategory of $\katez$. By $G$ we always denote a finite group. We will
consider faithful representations $\rho \colon G \longrightarrow \Aut_k k \psl t\psr.$
Let $A$ denote an object of $\kate$, and set $\Gamma_A = \Aut_A A\psl t \psr$. A \emph{deformation of $\rho$} to an
object $A$ of $\kate$ is a homomorphism $\tilde{\rho} \colon G
\to \Gamma_A$ such that the following diagram
commutes:
$$\xymatrix{& \Gamma_A
\ar[d] \\ G \ar[ur]^{\tilde{\rho}}
\ar[r]_{\rho} & \Gamma_k.}$$ Two deformations
are called \emph{equivalent} if they differ by conjugation by an element
of $\Gamma_{A,k}=\ker(\Gamma_A \to \Gamma_k)$.
The \emph{deformation functor} $$D_{\rho} \colon \kate \to \mathsf{{Sets}}$$
associates to $A$ the set of equivalence classes of deformations of
$\rho$ to $A$. We will often write $D$ for $D_{\rho}$ and also
denote by deformations the equivalence classes in which they lie. If we only consider lifts to rings $A$ in  $\kate$ of characteristic $p$, we arrive at the \emph{equicharacteristic deformation functor} that we denote by $D_\rho/p$.

The set $D(k[\varepsilon]/\varepsilon^2)$ has a structure of $k$-vector space and is called the tangent space to $D$. Using Schlessinger's criteria \cite{Schlessinger}, one
may easily prove  that for any $\rho$
the functor $D$ has a \emph{versal deformation ring} $R$ in $\katez$. This
means that there is a morphism of functors $h_R = \Hom(R,\cdot) \to
D$ which is smooth and induces an isomorphism on tangent spaces. 

The representation $\rho$ induces on the group $G$ a decreasing filtration by higher ramification groups
$ G \supseteq G_1 \supseteq G_2 \supseteq \dots $ with
$$G_{i}:=\{ \sigma \in G \ : \ \mbox{ord}_t(\rho(\sigma) t-t)>i \} \ \ (i \geq 1).$$
If $G_1=\{0\}$,  $\rho$ is \emph{tamely ramified}, if $G_2=\{0\}$,  $\rho$ is called \emph{weakly ramified}. For example, S.~Nakajima \cite{Nakajima} has shown that every action of a ramification group at a point of an ordinary curve is weakly ramified. Weak ramification is sometimes called ``Hasse-conductor one''.

Versal deformation rings for weakly ramified group actions were explicitly determined in \cite{BM} (cyclic $p$-group), \cite{CK} (general equicharacteristic case) and \cite{CM} (general case). For example, if $p \geq 3$, the action  of a cyclic $p$-group by $t \mapsto t/(1+t)$ has versal deformation ring $W(k)\psl \alpha \psr /\langle \psi(\alpha) \rangle$, where $\psi(\alpha)$ is a polynomial of degree $\frac{p-1}{2}$, and the versal deformation is given by $t \mapsto  ({t+\alpha})/({t+\alpha+1}). $
But is this deformation universal? 

\begin{introthm}
Let $\rho \colon G \to \Aut_k k\psl t \psr$ be a weakly ramified local representation of a finite group $G$. The pro-representability of $D_\rho$ only depends on the abstract type of the group $G$ and the characteristic $p$ of the ground field. 

More precisely, the functor $D_\rho$ is not  pro-representable if and only if $p=2$ and $G \in \{{\mathbf Z}/2, ({\mathbf Z}/2)^2\}$.\end{introthm}

\begin{introthm}
Let $\rho \colon G \to \Aut_k k\psl t \psr$ be a weakly ramified local representation of a finite group $G$. The pro-representability of $D_\rho/p$ only depends on the abstract type of the group $G$ and the characteristic $p$ of the ground field. 

More precisely, the equicharacteristic deformation functor $D_\rho/p$ is not  pro-representable if and only if $p=2$ and $G = \mathbf{Z}/2$.\end{introthm}

The question of universality is equivalent to that of injectivity of the map $h_R \rightarrow D$, i.e., to the following: suppose two morphisms $R \rightarrow A$ induce deformed representations to $A \in \kate$ that are conjugate by an element of $\Gamma_A$; then are these morphisms equal?

The proof is based on a kind of ``linearization'' technique, roughly as follows: call an element of $\Gamma_A$ a \emph{homography} if it is of the form $t \mapsto (at+b)/(ct+d)$; then ``when two homographies are conjugate by a power series, they are also conjugate by a homography''. Though this is false as it stands, the gist is right (see Lemma \ref{conjugation} for a correct statement, based on a more careful analysis of the ``Nottingham group'' over the category $\kate$).  Since all versal lifts of weakly ramified group actions are homographies, the observation allows one to reason in a much smaller space of conjugating objects and finish the proof: essentially because a homography is in general determined by its first three Taylor coefficients. 

The paper ends with a section that contrasts our method of proof with an abstract version of the general method employed in the existing literature, and we discuss a conjecture of Tim Dokchitser on non-pro-representable functors in the setting of a weak involution in characteristic two.

\begin{introrem}
Let $G$ denote a finite group acting on a  projective curve $X$ of genus $g \geq 2$. The ``global'' deformation functor $D_{X,G}$ of the pair $(X,G)$ admits a smooth morphism to the direct product of ``local'' deformation functors of the ramification groups at the completed local rings at ramification points \cite{BM}. The functor $D_{X,G}$ is \emph{always} pro-representable, for the simple reason that it has no so-called ``infinitesimal automorphisms'' (since $H^0(X,T_X)^G \subseteq H^0(X,T_X)=0$; see also \cite{Sernesi}, Section 2.6). This even holds in characteristic two, with the local deformation functors not necessarily pro-representable. 
\end{introrem}

\begin{section}{Nottingham and Borel groups over the category $\kate$}\label{matrices}

In this section, we set up the necessary technical lemmas that allow us to switch from inacessible calculations with power series to easier ones involving fractional linear transformations.

\begin{df} Let $A$ in  $\kate$, and denote by $\bar{\cdot} \colon A \rightarrow k$ the reduction map. We denote the group $\Aut_A A\psl t \psr$ by $\Gamma_A$. Define its subgroups $$\Gamma_A^i=\{\fhi \in \Gamma_A \mid \fhi(t) \equiv t \pmod{t^{i+1}} \}, \quad i \geq -1.$$ We obtain in this way a decreasing filtration $$\Gamma_A = \Gamma_A^{-1} \supseteq \Gamma_A^0 \supseteq \Gamma_A^1 \supseteq \dots.$$
An element $\gamma \in \Gamma_A$ is uniquely determined by the power series $\gamma(t)$, and by slight abuse of notation, we shall also write such an element of $\Gamma_A$ as the power series that represents its action on $t$. In terms of these power series, the group law corresponds to composition, $(\varphi \cdot \psi)(t)=\psi(\varphi(t))$.
\end{df} 

\begin{rem}
The group of ``wild automorphisms'' $\Gamma^0_A$ has been called the Nottingham group in group theory, especially for $A$ a finite field , cf.\ \cite{Nottingham} chapter 6 \& 10. Here, we study the ``Nottingham group over the category $\kate$''. We need to be careful with the ``linear algebra'' over rings $A $ in $\kate$, since they can have nilpotents, etc. Also note that for all $A \neq k$ in $\kate$, there exist $\gamma \in \Gamma_A$ such that $\gamma(t)$ has non-zero constant term. 
\end{rem}

\begin{lemma}\label{easylem} For $i \geq 1$, we have \begin{enumerate} \item[\textup{(i)}] For a $\varphi \in \Gamma_A$, we have $\varphi(t)=a_0+a_1 t + \dots$ with $a_0 \in \mathfrak{m}_A$ and $a_1 \in A^*$;
\item[\textup{(ii)}] For any $g(t) \in A\psl t \psr$ such that $g(t)=a_0+a_1 t + \dots$ and $a_0 \in \mathfrak{m}_A$, $a_1 \in A^*$, there exists a unique $\varphi \in \Gamma_A$ such that $\varphi(t)=g(t)$;
\item[\textup{(iii)}] $\psi \in \Gamma_A$ belongs to the left coset $\Gamma_A^i \gamma$ exactly if $\psi(t) \equiv \gamma(t) \pmod{t^{i+1}}$;
\item[\textup{(iv)}] the subgroups $\Gamma_A^{i}$ are normal in $\Gamma_A^0$. \end{enumerate} \end{lemma} 

\begin{proof}
Statements (i) and (ii) are easy. Concerning the second statement, for $\gamma = \sum b_j t^j \in \Gamma_A$, a series $\psi=\sum a_j t^j$ belongs to the left coset $\Gamma_A^i \gamma$ exactly if there exists $\gamma_i = t + \sum_{j \geq i+1} c_j t^j$ with $\gamma_i \cdot \gamma = \psi$, which is is equivalent to an infinite linear system of the form
$a_0=b_0,\, a_1=b_1,\, \dots ,\, a_i=b_i$ and for $j>i$, $$a_j = b_1 c_{j} +\mbox{ (an algebraic expression in $c_k$ for $k<j$ and $b_k$)}.$$ Now the important thing is that $b_1$ is invertible (part (i)), so one can solve iteratively for $c_j$. 

For  (iii), we determine the right cosets $\gamma \Gamma_A^i$ for $\gamma \in \Gamma_A^0$ in a similar way: since then $b_0=0$, we find an infinite linear system of the form $a_0=0,\,  a_1=b_1,\, \dots ,\, a_i=b_i$ and for $j>i$, $$ a_j = c_j b_1^{i+1} +\mbox{ (an algebraic expression in $c_k$ for $k<j$ and $b_k$)}.$$ Again, $b_1$ is invertible (part (i)), so one can solve iteratively for $c_j$. Since left and right cosets coincide, $\Gamma_A^i$ is normal in $\Gamma_A^0$. 
\end{proof}

\begin{rem}
The groups $\Gamma_A^{i}$ ($i \geq 0$) are in general \emph{not} normal in $\Gamma_A$, as can be seen from constructing the right cosets of  a $\gamma=b_0+b_1 t + \dots \in \Gamma_A$ with $b_0 \neq 0$ as in the above proof. Another proof of (iii) goes as follows: $\Gamma_A^i$ is the kernel of the morphism from $\Gamma_A^0$ to the power series truncated at $t^i$. For $A$ a finite field, (ii) and (iii) are also proven on page 207--208 of \cite{Nottingham}. 
\end{rem}

\begin{df} We denote by $\PB(A)$ the subgroup of the group $\PGL_2(A)$ given by $$\PB(A)=\left\{\Big(\begin{matrix} a&b \\c&d 
\end{matrix}\Big) \in \PGL_2(A) \mid a,c,d \in A, \; b \in \mathfrak{m}_A \right\}.$$\end{df}

\begin{rem} For $A=k$, $\PB(k)$ is precisely the standard (projective) Borel subgroup $\PB(k)=B(k)/k^*$ of the group $\PGL_2(k)$ and thus the elements of $\PB(A)$ can be regarded as deformations of matrices in $\PB(k)$. Note that $\PB$ is a group functor on the category of local rings with residue field $k$.
\end{rem}

\begin{df}\label{przeobrazenie}
The reduction map $\bar{\cdot} \, : \, A \rightarrow k$ induces reduction maps $\Gamma_A \rightarrow \Gamma_k$ and  $\PB(A) \rightarrow \PB(k)$ that we still denote by a bar. Their respective kernels are denoted by $\Gamma_{A,k}$ and $\PB_{A,k}$.
\end{df}

\begin{se} To any element $\big( \begin{smallmatrix} a&b \\ c&d \end{smallmatrix} \big) \in \PB(A)$ we associate the element of $\Gamma_A$ given by its Taylor expansion \begin{equation}\label{fractransfo} t \mapsto \frac{at+b}{ct+d} = \frac{b}{d} + \frac{ad-bc}{d^2}t-\frac{c(ad-bc)}{d^3}t^2 + \ldots.\end{equation} In this way, we can interpret $\PB(A)$ as a subgroup of $\Gamma_A$. We will do so without further mention. This interpretation is functorial in $A$. \end{se}

\begin{lemma} \label{rozwielitka} Any element $\gamma \in \Gamma_A$ has a unique decomposition of the form $$\gamma = \gamma_2 \fhi \quad \text{with} \quad \gamma_2 \in \Gamma_A^2, \fhi \in \PB(A).$$
\end{lemma}
\begin{proof} By Lemma \ref{easylem} (iii), we can represent the elements of $\Gamma_A^2 \backslash \Gamma_A$ by power series modulo $t^3$.
The surjective set theoretical map $f : \PB(A) \rightarrow \Gamma^2_A \backslash \Gamma_A$ induced by (\ref{fractransfo}) has inverse $$u_0 + u_1 t + u_2 t^2 \ \mathrm{ mod }\ t^3 \mapsto \Big(\begin{matrix} u_1+u_0u_2u_1^{-1}&u_0 \\ u_2u_1^{-1}&1
\end{matrix}\Big).$$ Note that the right hand side belongs to $\PB(A)$ since $u_0 \in \mathfrak{m}_A$ and its determinant is $u_1 \in A^*$, cf.\ Lemma \ref{easylem}(ii). 

For general $\gamma \in \Gamma_A$, set $\varphi:=f^{-1}(\Gamma^2_A \gamma)$; then $\gamma_2:=\gamma \varphi^{-1} \in \Gamma_A^2$, and the uniqueness of the decomposition follows from the bijectivity of $f$. 
\end{proof}

From this, we deduce our first main property of conjugation of elements in $\PB(A)$ by an element of $\Gamma_A$ (recall the definition of $\Gamma_{A,k}$ and $\mathfrak{B}_{A,k}$ from \ref{przeobrazenie}):

\begin{prop}\label{conjugation} Let $\fhi,\psi \colon G \to \PB(A)$ be two morphisms of groups. Assume that $\psi$ can be conjugated into $\Gamma_A^0$ by an element  in $\PB_{A,k}$. 
Then if $\fhi$ and $\psi$ are conjugate by an element of $\Gamma_{A,k}$, then they are also conjugate by an element of $\PB_{A,k}$. 
\end{prop} 
\begin{proof} We assume there exists $\tau \in \PB_{A,k}$ and $\eta \in \Gamma_{A,k}$ with $\tau \psi(g) \tau^{-1} \in \Gamma^0_A$ and  $\psi(g) = \eta \fhi(g) \eta^{-1}$ for all  $g\in G.$  
Hence we know that $\tau \psi(g) \tau^{-1} = \tau\eta \fhi(g) (\tau\eta)^{-1}$ lies in $\Gamma^0_A$ for any $g\in G$ and that $\overline{\tau \eta}=\id$. By Lemma \ref{rozwielitka} we can write $\tau \eta = \xi \theta$ with $\xi \in \Gamma^2_A$ and $\theta \in \PB(A)$. Furthermore, by the uniqueness of such a decomposition over $k$ we have $\overline{\xi}=\overline{\theta}=\id$. Now $$\xi \cdot \theta \fhi(g) \theta^{-1} = \tau \psi(g)\tau^{-1} \cdot \xi.$$ Since $\Gamma_A^2$ is a normal subgroup of $\Gamma_A^0$ (cf.\ Lemma \ref{easylem}(iv)), we have $$\tau \psi(g) \tau^{-1} \cdot \xi = \xi'_g \cdot \tau \psi(g)\tau^{-1}$$ for some $\xi'_g \in \Gamma^2_A$ (which \emph{a priori} might depend on $g$). Thus $$\xi \cdot \theta \fhi(g) \theta^{-1} = \xi'_g \cdot \tau \psi(g)\tau^{-1}$$ are two decompositions into a product of an element of $\Gamma^2_A$ and $\PB(A)$. Again by Lemma \ref{rozwielitka} they coincide, i.e., $\xi = \xi'_g$ and $\tau \psi(g) \tau^{-1} = \theta \fhi(g) \theta^{-1}$. Thus $$ \psi(g) = \tau^{-1}\theta \fhi(g) (\tau^{-1}\theta)^{-1}$$ for all $g\in G$ and $\tau^{-1}\theta \in \PB(A)$, $\overline{\tau^{-1}\theta}=\id$.
\end{proof}

\begin{rem} \label{bad}
In the proposition, the condition that one of the representation can be conjugated into $\Gamma_A^0$ cannot be left out (in the proof, this is reflected in the use of the fact that $\Gamma_A^2$ is normal in $\Gamma_A^0$). 

Indeed, let $k$ be a field of characteristic $\kar k \neq 2,3$ and set $A=k[\varepsilon]/\varepsilon^3$. Then $$\psi \colon t \mapsto t+\varepsilon \mbox{ and }\fhi \colon t \mapsto \frac{t+\varepsilon}{-3\varepsilon^2t+1}=\varepsilon + t + 3 \varepsilon^2 t^2$$ are conjugate by $$\tau : t \mapsto t + \varepsilon t^3 \in \Gamma_A^0$$ since $\psi \tau = \tau \fhi$, but are not conjugate by an element $\gamma$ of $\PB(A)$. The latter is proven by direct matrix calculation, along the following lines. Suppose we have a matrix representation $\gamma \psi = \lambda \fhi \gamma$ for $\lambda \in A^*$. Taking traces, we find $\lambda=1$, and then we arrive at a system of equations for the entries of $\gamma$ that only has a solution with all those entries in $\maks_A$, so such $\gamma$ are not invertible.
And indeed, neither $\psi$ nor $\fhi$ can be conjugated into $\Gamma_A^0$ --- we leave out the standard matrix computation that proves this.
\end{rem}

The next proposition shows that certain group actions can be conjugated into a $\Gamma^0$-type group, but only over an extension of the ground ring $A$: 

\begin{prop}\label{wiet} Let $G$ be a cyclic $p$-group and $\fhi \colon G \to \PB(A)$ a group homomorphism. There exists an extension $A \subseteq A'$ in $\kate$ such that $\fhi$ can be conjugated by an element of $\PB_{A',k}$ into $\Gamma_{A'}^0$.
\end{prop}
\begin{proof}  Let $g$ denote a generator of $G$ and put $$A'=A[\mathfrak{z}]/(\mathfrak{z}^2-a\mathfrak{z}-b)$$ with some $a,b \in \mathfrak{m}_A$ to be determined. The ring $A'$ is a free $A$-algebra with basis $\{1,\mathfrak{z}\}$. Since $\mathfrak{z}^2 A' \subseteq \mathfrak{m}_A A'$, the ideal $\mathfrak{m}_A \cdot 1 + A \cdot \mathfrak{z}$ is nilpotent. Hence, it is the unique maximal ideal and the ring $A'$ is an object of $\kate$, in particular, the residue field is $k$. Write $\fhi(g) = \big( \begin{smallmatrix} u&v \\ w&z \end{smallmatrix} \big)$ and put $\gamma(t)=t-\mathfrak{z}$. Then $\gamma \fhi(g) \gamma^{-1} \in \Gamma_{A'}^0$ if and only if $$\frac{u \mathfrak{z} + v}{w\mathfrak{z}+z}=\mathfrak{z}.$$ This means that we need to have $wa=u-z$ and $wb=v$. The order of $\bar{\fhi}(g)=\big( \begin{smallmatrix} \bar{u}&0 \\ \bar{w}&\bar{z} \end{smallmatrix} \big)$ is a power of $p$, and is different from one. Hence we have $\bar{u}=\bar{z}$ and $\bar{w} \neq 0$. Thus $w$ is invertible and $a=w^{-1}(u-z)$ and $b=w^{-1}v$ give the unique solution to our equations. 
\end{proof}

\begin{rem}
The extension $A \subseteq A'$ cannot be avoided in general. For example, with $\kar k = 2$ and $A=k[\varepsilon]/\varepsilon^4$, we have that$$ \Big(\begin{matrix} 1 &\varepsilon \\1&1 
\end{matrix}\Big) \mbox{  and  } \Big(\begin{matrix} 1&\varepsilon+\varepsilon^3 \\1&1 
\end{matrix}\Big) $$ are not conjugate by $\PB_{A,k}$, but they are by $\PB_{A',k}$ for $A'=A[\sqrt{1+\varepsilon^3}]$. This can be seen by an easy calculation, or by appealing to Proposition \ref{psychopatologia} \emph{infra}. 
\end{rem}

\begin{remark} The proof shows more: namely, that we can choose $A'=A[\mathfrak{z}]/(\mathfrak{z}^2-a\mathfrak{z}-b)$ with $a,b \in \mathfrak{m}_A$ and $\gamma$ to be just the translation $\gamma(t)=t-\mathfrak{z}$.
\end{remark}

We wish to extend this proposition to the case of certain abelian $p$-groups, for which we need the following:

\begin{df}
Two matrices $\psi, \fhi$ in $\PGL(2,A)$ are called \emph{affinely dependent} if there exist $a \in A^*$ and $b \in A$ such that $\psi = a \fhi + b \cdot \mathrm{id}$. The relation of affine dependence is an equivalence relation. 
\end{df}

\begin{lem} \label{commute} Assume $\kar k \neq 2$.  Then any two commuting matrices $\fhi,\psi$ in $\PB(A)$ such that their images in $\PB(k)$ are not both diagonal are affinely dependent.
\end{lem}
\begin{proof}
Let $m$ and $n$ be matrices inducing $\fhi$ and $\psi$, with $\bar{m}$ not diagonal. Then there exists $\lambda \in A^*$ such that $mn=\lambda nm$. Taking determinants we get $\lambda^2 =1$. Now $(\lambda+1)-(\lambda-1)=2$, and since $\kar k \neq 2$, we conclude that one of $\lambda \pm 1$ is invertible. Hence we conclude from $(\lambda+1)(\lambda-1)=0$ that $\lambda = \pm 1$.

Write $m=\big( \begin{smallmatrix} a_1&b_1 \\ c_1&d_1 \end{smallmatrix} \big)$ and $n=\big( \begin{smallmatrix} a_2&b_2 \\ c_2&d_2 \end{smallmatrix} \big)$.
This gives $$\Big(\begin{matrix} a_1a_2+b_1c_2&a_1b_2+b_1d_2 \\ c_1a_2+d_1c_2&c_1b_2+d_1d_2 \end{matrix}\Big)=\lambda \Big(\begin{matrix} a_1a_2+c_1b_2&b_1a_2+d_1b_2 \\ a_1c_2+c_1d_2& b_1c_2+d_1 d_2 \end{matrix}\Big).$$ Since $\bar{a_1}, \bar{a_2} \in k^*$ and $\bar{b_1}=\bar{b_2}=0$, by reduction to $k$ one gets $\bar{\lambda}=1$ from looking at the left most top entry. Hence $\lambda =1$. It then follows that $b_1c_2=c_1b_2$ and $(a_1-d_1)c_2=c_1(a_2-d_2)$.  As the reduction of $m$ is not diagonal, $\bar{c_1} \neq 0$, so $c_1$ is invertible, and this shows that $$n=am+b \cdot \mathrm{id} \mbox{ with } a=\frac{c_2}{c_1} \mbox{ and } b=\frac{c_1d_2-d_1c_2}{c_1}.$$
\end{proof}

\begin{rem} \label{quol}
The claim of the lemma is false if $\kar k=2$, for example, in $A=k[\varepsilon]/\varepsilon^2$,  set $$ \psi := \Big(\begin{matrix} 1+\varepsilon&\varepsilon \\0&1 
\end{matrix}\Big) \mbox{  and  } \fhi:= \Big(\begin{matrix} 1&\varepsilon \\1&1 
\end{matrix}\Big), $$  then  $\psi \fhi = \lambda \fhi \psi$, for $\lambda = 1 + \varepsilon$, and $\bar{\fhi}$ is not diagonal, but nevertheless $\psi \neq a \fhi + b \cdot \mathrm{id}$ for any $a \in A^*, b \in A$. And indeed, $\lambda^2=1$ but $\lambda \neq \pm 1$. 
\end{rem}

Now comes the desired extension of the previous proposition:

\begin{prop}\label{elab} Let $G$ be an abelian $p$-group with $p \neq 2$ and $\fhi \colon G \to \PB(A)$ a group homomorphism. There exists an extension $A \subseteq A'$ in $\kate$ such that $\fhi$ can be conjugated by an element of $\PB_{A',k}$ into $\Gamma_{A'}^0$.
\end{prop}

\begin{proof} If $G$ is an abelian $p$-group and $p \neq 2$, note that the only diagonal matrix in the reduction of $\fhi(G)$ is the identity matrix. Since the image of $\fhi$ lies in $\PB(A)$, Lemma \ref{commute} applies, so all the elements of the image are affinely dependent. Hence the conjugation from Proposition \ref{wiet} applies simultaneously to all elements of $G$.
\end{proof}

We will also need the following special form in characteristic two, where we have to make the property in Lemma \ref{commute} an extra hypothesis:

\begin{prop}\label{kar2} Let $G$ be an abelian $2$-group and $\fhi \colon G \to \PB(A)$ a group homomorphism such that all elements in the image of a set of generators for $G$ under $\fhi$ are pairwise affinely dependent. Then there exists an extension $A \subseteq A'$ in $\kate$ such that $\fhi$ can be conjugated by an element of $\PB_{A',k}$ into $\Gamma_{A'}^0$. 
\end{prop}

\begin{proof}
This is literally the same as that of Proposition \ref{elab}, except that the property of affine dependence is now an assumption.
\end{proof}

We end this section by a proposition about conjugacy in $\PB(A)$.

\begin{df}
Let $m$ and $n$ denote two matrices in $\PGL_2(A)$. We define an equivalence relation $m \approx n$ by the existence of a residually trivial conjugacy between $m$ and $n$, i.e., by $\chi m \chi^{-1} = n$ for some $\chi \in \PGL_2(A)$ with $\bar{\chi}=\id$. 
\end{df}

\begin{prop}\label{psychopatologia}  Let $m$ and $n$ be two matrices in $\PB(A)$ such that $\bar{m}$ and $\bar{n}$ are not diagonal. Then the following conditions are equivalent: \begin{enumerate} \item[\textup{(i)}] $m \approx n$; \item[\textup{(ii)}] There exists representatives $\tilde{m}$ and $\tilde{n}$ in $\GL_2(A)$ for $m$ and $n$ such that $$ \tr \tilde{m} = \tr \tilde{n}, \ \det \tilde{m} = \det \tilde{n} \mbox{ and } \bar{\tilde{m}}=\bar{\tilde{n}}. $$
\end{enumerate}
Furthermore, if $\tr m \notin \maks_A$, then \textup{(ii)} is equivalent to
\begin{enumerate}
\item[\textup{(ii')}] $\displaystyle \frac{(\tr m)^2}{\det m}  = \frac{(\tr n)^2}{\det n}$ and $\bar{m}=\bar{n}$. 
\end{enumerate}
\end{prop}
\begin{proof} It is clear that (i) implies (ii). A guiding principle for the proof in the other direction is the following: if $A$ were a field, we have two matrices with the same characteristic polynomial, hence their Jordan Normal Forms are equal, so they are conjugate. For $A \in \kate$, we cannot use this argument, hence we replace it by an explicit reduction of matrices.

For ease of notation, we use the letters $m,n$ for representatives in $\GL_2(A)$ of the given matrices, and we use $\approx$ for matrices in $\GL_2(A)$ to mean conjugacy by a matrix whose reduction is trivial. Set $m=\big( \begin{smallmatrix} a_1&b_1\\ c_1&d_1 \end{smallmatrix} \big)$ and $n=\big( \begin{smallmatrix} a_2&b_2\\ c_2&d_2 \end{smallmatrix} \big)$.

We observe the following useful identity \begin{equation}\label{gorzalka}\Big( \begin{matrix} 1&\alpha \\ 0&1 \end{matrix} \Big) \Big( \begin{matrix} a&b \\ c&d \end{matrix} \Big) \Big( \begin{matrix} 1&-\alpha \\ 0&1 \end{matrix} \Big) = \Big( \begin{matrix} a+\alpha c&b + \alpha (d-a) - \alpha^2 c\\ c&d-\alpha c \end{matrix} \Big).\end{equation} 

With $\alpha:= c_1^{-1} (d_1-d_2)$, we find from (\ref{gorzalka}) that $$ m \approx m':= \Big( \begin{matrix} a_2&b_1' \\ c_1&d_2 \end{matrix} \Big),$$
with some $b_1'$ such that $\bar{b_1'}=0$. Indeed, observe that  $\alpha \in \maks_A$, so that $\bar{m'}=\bar{m}$. Thus, condition (ii)  continues to hold also for $(m',n)$. 

Now write $\mu=c_2/c_1 \in A^*$. Then $\overline{\mu}=1$ since $\bar{m}=\bar{n}$. We have $$\Big( \begin{matrix} \mu^{-1}&0\\ 0&1 \end{matrix} \Big) m' \Big( \begin{matrix} \mu^{-1} &0 \\ 0&1 \end{matrix} \Big)^{-1} = m'':= \Big( \begin{matrix} a_2 & b_1'' \\ c_2 &d_2 \end{matrix} \Big)$$
for some $b_1''$ and (ii) still holds for $(m'',n)$. Since $c_2$ is invertible, the determinant condition implies that $b_1''=b_2$, and this proves that $m \approx n$.  

Finally, we consider condition (ii').  In general, (ii) implies (ii'). Conversely, if $\tr m \notin \maks$, then also $\tr n \notin \maks$, and with $\lambda:= \tr m / \tr n \in A^*$, the pair of representatives $(m, \lambda n)$    
satisfies (ii), since $\bar{\lambda}=1$.
\end{proof}

\begin{rem}
Conditions (ii) and (ii') in the proposition are not equivalent if $\tr m \in \maks_A$, for example, with $A=k[\varepsilon]/\varepsilon^2$, set $$m=\Big( \begin{matrix} 1+\varepsilon & 0\\ 1& -1 \end{matrix} \Big) \mbox{ and } n=\Big( \begin{matrix} 1&0\\ 1&-1 \end{matrix} \Big),$$ then (ii') holds, but (ii) doesn't.

The equivalence of (i) and (ii) can also fail if the reduced matrices are diagonal, for example, the identity matrix and the matrix $\big( \begin{smallmatrix} 1 & \varepsilon\\ \varepsilon & 1 \end{smallmatrix} \big)$ over $k[\varepsilon]/\varepsilon^2$ 
clearly satisfy (ii) but are not conjugate.

\end{rem}

\end{section}

\begin{section}{Deformation of the action of a cyclic group of order $p$}

\begin{se} In this section we study the example case where $G={\mathbf Z}/p$ for $\kar k = p \neq 2$. Let $g$ denote a generator of $G$ and suppose $\rho: G \rightarrow \Gamma_k$ is weakly ramified. Since the Hasse conductor is one, by Artin-Schreier theory, we can normalize the action of $G$ on $t$ to be of the form $t \mapsto t/(1+t)$: the corresponding field extension $k(\!(t)\!)/k(\!(x)\!)$ is given by $(1/t)^p-1/t=1/x$ with Galois group generated by $1/t \mapsto 1/t+1=t/(1+t)$. We recall the form of the versal deformation from \cite{BM}: the versal deformation ring of $D_\rho$ is given by $R=W(k)\psl\alpha\psr /\langle \psi(\alpha) \rangle$, where $$\psi(\alpha)=\sum_{l=0}^{\frac{p-1}{2}} \binom{p-1-l}{l} (-1)^l (\alpha +4)^{\frac{p-1}{2}-l}$$ and the versal deformation is given by $$\tilde{\rho}(g)(t) = \frac{t+\alpha}{t+\alpha+1}.$$ 
\end{se}

\begin{prop} \label{Z/p}
For $p \neq 2$, $G$ a group of order $p$ generated by $g$, and the action $\rho : G \rightarrow \Gamma_k$ with $\kar k = p$ given by $\rho(g)=t/(1+t)$, the local deformation functor $D_\rho$ is pro-representable.
\end{prop}

\begin{proof} Let $R$ denote the above versal deformation ring. 
Assume that for $A$ in $\kate$, the map $\Phi_A \colon h_R(A) \to D_{\rho}(A)$ is not injective. Choose $\fhi_1, \fhi_2 \in h_{R}(A)$ with the same image in $D_{\rho}(A)$. Write $\alpha_1=\varphi_1(\alpha)$, $\alpha_2=\varphi_2(\alpha)$, and let $$m_i := \varphi_i^* \circ \tilde{\rho}(g) =   \Big(\begin{matrix} 1&\alpha_i \\ 1&\alpha_i+1
\end{matrix}\Big). $$ 
The assumption is that $m_1$ and $m_2$ are conjugate in $\Gamma_A$. 
Applying Proposition \ref{wiet} to $\varphi_2^* \circ \tilde{\rho}$, we know that $m_2$ can be conjugated by $\PB_{A',k}$ into $\Gamma_{A'}^0$ for some extension $A \subseteq A'$ in $\kate$. Then, since $m_i \in \PB(A')$, by Lemma \ref{conjugation}, we know that $m_1$ and $m_2$ are conjugate by an element of $\PB_{A',k}$, in particular, $m_1 \approx m_2$. Now since $\tr m_2 = 2 + \alpha_2 \notin \maks_{A'}$ (as $p \neq 2$), we find from condition (ii') of Proposition \ref{psychopatologia} that $$0 = (\tr m)^2-(\tr n)^2 = (\alpha_1+2)^2-(\alpha_2+2)^2=(\alpha_1+\alpha_2+4)(\alpha_1-\alpha_2).$$ Since residually $\bar{\alpha_1}=\bar{\alpha_2}=0$ but $\bar{4} \neq 0$,  $\alpha_1 + \alpha_2+4$ is a unit, and thus the above equation implies $\alpha_1=\alpha_2$, so $\fhi_1=\fhi_2$, a contradiction.  
\end{proof}

\end{section}

\begin{section}{Non-pro-representable deformation functors}

We postpone the treatment of the other pro-representable cases of deformation functors of weakly ramified local actions to the next section, since the arguments are just a technical enhancement of those in the previous section. We first concentrate on the anomalous cases where $D_\rho$ is not pro-representable.

\begin{se} \label{two}
Let $G={\mathbf Z}/2$, $\kar k = 2$; let $g$ denote a generator of $G$ and suppose $\rho: G \rightarrow \Gamma_k$ is weakly ramified. By Artin-Schreier theory, we can normalize the action of $G$ on $t$ to be of the form $t \mapsto t/(1+t)$. We recall the form of the versal deformation from \cite{BM}: the versal deformation ring of $D_\rho$ is given by $R=W(k)\psl\alpha\psr$, and the versal deformation is given by $$\tilde{\rho}(g)(t) = \frac{t+\alpha}{t-1}.$$ 
\end{se}

\begin{prop}
For $G$ a group of order $2$ generated by $g$ and the action $\rho : G \rightarrow \Gamma_k$ ($\kar k = 2$) given by $\rho(g)=t/(1+t)$, the local deformation functors $D_\rho$ and $D_\rho/2$ are not pro-representable.
\end{prop}

\begin{proof} Let $R$ denote the above versal deformation ring. Then $R/2$ is the versal deformation ring for $D_\rho/2$. 
Take $A=k[\varepsilon]/\varepsilon^3$ and let $\fhi_1, \fhi_2 \colon W(k)\psl\alpha\psr \to A$ be homomorphisms given by $\alpha \mapsto \varepsilon$ and $\alpha \mapsto \varepsilon + \varepsilon^2$ respectively. The map $h_{R/2}(A) \to D(A)$ maps $\fhi_1$ and $\fhi_2$ to deformations $\rho_1, \rho_2 \in D(A)$ given by $$\rho_1(g)(t) = \frac{t+\varepsilon}{t-1} \mbox{ and }\rho_2(g)(t) = \frac{t+\varepsilon+\varepsilon^2}{t-1},$$ respectively. However, the equality $$ (1+\varepsilon) \Big(\begin{matrix} 1+\varepsilon&\varepsilon \\0&1 
\end{matrix}\Big) \Big(\begin{matrix} 1&\varepsilon \\1&-1 
\end{matrix}\Big) \Big(\begin{matrix} 1+\varepsilon&\varepsilon \\0&1 
\end{matrix}\Big)^{-1} = \Big(\begin{matrix} 1&\varepsilon + \varepsilon^2 \\1&-1 
\end{matrix}\Big)$$ shows that those two deformations are equivalent. Thus, the map $h_{R/2}(A) \to D(A)$ is not injective, and neither the ring $R$ nor $R/2$ is universal. 
\end{proof}

\begin{se} We proceed with the next non-pro-representable case:
let $G=(\mathbf{Z}/2)^2$. This case is special, since though the mixed-characteristic deformation functor is not pro-re\-pre\-sen\-ta\-ble, the equicharacteristic one is.

 We can normalize any weakly ramified $\rho \colon G \to \Gamma_k$ to be given by $$\rho(1)(t)=\frac{t}{t+1},  \ \ \rho(u)(t)=\frac{t}{ut+1}$$ for $\langle 1, u \rangle \subseteq k$ an ${\mathbf F}_2$-vectorspace of dimension 2 (cf. \cite{CKgerman}). The versal deformation ring of $\rho$ is given by $R=W(k)\psl\alpha\psr$ and the versal deformation is given by $$\tilde{\rho}(1)(t) = \frac{t+\alpha}{t-1}, \quad \tilde{\rho}(u)(t) = \frac{t-(\alpha\tilde{u}+2)}{\tilde{u}t-1},$$ with $\tilde{u}$ any lift of $u$ to $W(k)$ (cf. \cite{CM}, proof of Proposition 3.8.i).  
\end{se}

\begin{prop}
For $G = ({\mathbf Z}/2)^2$ and a weakly ramified action $\rho : G \rightarrow \Gamma_k$  with $\kar k = 2$, the local deformation functor $D_\rho$ is not pro-representable. 
\end{prop}

\begin{proof} Let $R$ denote the above versal deformation ring. 
Take $A=W(k)/16$ and let $$\fhi_1, \fhi_2 \colon W(k)\psl\alpha\psr \to A$$ be homomorphisms given by $\alpha \mapsto -2$ and $\alpha \mapsto 6$ respectively. Put $$ \gamma = \Big(\begin{matrix} 5-4\tilde{u}(1-\tilde{u})&-4-4\tilde{u}(1-\tilde{u}) \\2\tilde{u}(1-\tilde{u})&1 
\end{matrix}\Big).$$ One checks then that the following equations hold

$$ \left\{ \begin{array}{l}  5\gamma  \Big(\begin{matrix} 1&-2 \\1&-1 
\end{matrix}\Big) \gamma^{-1} = \Big(\begin{matrix} 1&6 \\1&-1 
\end{matrix}\Big), \\ (5-4(\tilde{u}-1)^2(2\tilde{u}+1))\gamma \Big(\begin{matrix} 1&2\tilde{u}-2 \\\tilde{u}&-1 
\end{matrix}\Big)\gamma^{-1} = \Big(\begin{matrix} 1&-6\tilde{u}-2 \\\tilde{u}&-1 
\end{matrix}\Big).\end{array} \right. $$
Hence the deformations given by $\fhi_1$and $\fhi_2$ are equivalent, so  the map $h_{R}(A) \to D_\rho(A)$ is not injective.
\end{proof}

However, we have the following:

\begin{prop} \label{smut}
For $G = ({\mathbf Z}/2)^2$ and a weakly ramified action $\rho : G \rightarrow \Gamma_k$  with $\kar k = 2$, the local equicharacteristic deformation functor $D_\rho/2$ \emph{is} pro-representable.
\end{prop}

\begin{proof} The ring $R/2$ is versal for $D_\rho/2$. Assume that $A \in \kate$ is of characteristic $2$ such that the map $\Phi_A \colon h_{R/2}(A) \to D_{\rho}(A)$ is not injective. Choose $\fhi_1, \fhi_2 \in h_{R/2}(A)$ with the same image in $D_{\rho}(A)$. Write $\alpha_1=\varphi_1(\alpha)$, $\alpha_2=\varphi_2(\alpha)$, and let $$m_i := \varphi_i^* \circ \tilde{\rho}(1) =   \Big(\begin{matrix} 1&\alpha_i \\ 1&1
\end{matrix}\Big); \ \ n_i := \varphi_i^* \circ \tilde{\rho}(1) =   \Big(\begin{matrix} 1&u \alpha_i \\ u &1
\end{matrix}\Big). $$ 
The assumption is that $m_1$ and $n_1$ are simultaneously  conjugate to $m_2$ and $n_2$ in $\Gamma_A$. 
Now $\varphi_2^* \circ \tilde{\rho}$, i.e., $m_2$ and $n_2$ can be conjugated into $\Gamma^0_A$: for this, we can use Lemma \ref{kar2}, since $m_i$ and $n_i$ ($i=1,2$) are affinely dependent via $n_i = u m_i + (u+1)\cdot \mathrm{id}$. By Proposition \ref{conjugation}, there exists a matrix $\gamma=\big( \begin{smallmatrix} a & b \\ c & d \end{smallmatrix} \big)$ with $\bar{\gamma}=\mathrm{id}$ such that  $$ \left\{ \begin{array}{l} m_2 \gamma = \lambda \gamma m_1  \\ n_2 \gamma= \mu \gamma n_1  \end{array} \right. $$ holds for some $\lambda,\mu \in A^*$. 
The equations of affine dependency imply $$u(\lambda-\mu)m_1=(u+1)(\mu-1) \cdot \mathrm{id}.$$ The left bottom entry of this matrix equation gives $\lambda=\mu$ and then, the equation implies $\lambda=\mu=1$. 
Hence the bottom row of the matrix equality $m_2 \gamma = \gamma m_1$  implies $a=d$ and $b=\alpha_1 c$ and then its top right entry gives $a(\alpha_1-\alpha_2)=0$. Since $a$ is invertible, we get $\alpha_1=\alpha_2$, and thus $\varphi_1=\varphi_2$. 
\end{proof}

\begin{rem}
What fails in general characteristic\footnote{Recall that the characteristic of a \emph{ring} $A$ is the unique nonnegative generator of the kernel of the natural map $\Z \to A$, cf. \cite{Bo} \S8, no. 8.}  (e.g., $\kar A=4$) is exactly the affine dependence of $\tilde{\rho}(1)$ and $\tilde{\rho}(u)$ that is so crucial in the above proof. 
\end{rem}

\begin{rem}\label{qual}
The proof shows that a matrix $\gamma$ with $\bar{\gamma}=\mathrm{id}$ that projectively commutes with $m_1$ \emph{and} $n_1$ is affinely dependent on $m_1$. Note that Remark \ref{quol} shows that just commutation with $m_1$ alone doesn't imply this.
\end{rem}

\end{section}

\section{Further pro-representable cases}

This section will consist of an enumeration (up to normalization, cf.\ \cite{CKgerman}) of all further possible weakly ramified actions, and proofs of the universality of their versal deformation rings. The fact that this list is complete and the versal rings are as indicated is the main contents of \cite{CM}. We will also take the occasion to clarify some points in that reference and in \cite{CK}. 

\begin{se}[$G={\mathbf Z}/n \mbox{ with } (n,p)=1 \mbox{ or } p=2, G=A_4$] In these cases, the versal deformation ring is $R=W(k)$ (the deformation problem is ``rigid''), hence $h_R(A)$ always consists of precisely one element, so $h_R(A) \rightarrow D_\rho(A)$ is necessarily injective and $D_\rho$ is pro-representable. \qed
\end{se}

\begin{rem}
There is a misprint in the versal deformation of $A_4$ on the bottom of page 251 in \cite{CM}. The correct unique lift to characteristic zero of $A_4$ is given by the following elements of $\PGL(2,W(k))$:
$$  m = \Big( \begin{matrix} 1 & 2 \\ 1 & -1 \end{matrix} \Big), \ \ n = \Big( \begin{matrix} 1 &-2j-2 \\ j & -1 \end{matrix} \Big), \ \ g = \Big( \begin{matrix} 1 & 0 \\ 0 & j \end{matrix} \Big), $$
with $j^2+j+1=0$. Then $m,n$ commute and are of order two, $g$ is of order three and $gmg^{-1}$ is equivalent to $n$. Also note that $j$ exists in $W(k)$ if $A_4$ is to have a residual representation over $k$, so it doesn't need to be adjoined, as is mistakenly done in the table on page 253 of \cite{CM}. 
\end{rem}

\begin{se}[$p \geq 3, G=D_p$] The versal deformation ring is the same as for a cyclic $p$-subgroup, and the argument of Proposition \ref{Z/p} applies literally. \qed
\end{se}

\begin{se}[$p \geq 3, G=({\mathbf Z}/p)^t \mbox{ or } G=({\mathbf Z}/p)^t \rtimes {\mathbf Z}/2 \mbox{ with } t \geq 2$] We have the following normalization of the residual weakly ramified representation: $(\mathbf{Z}/p)^t$ can be seen as a vector space $V$ of dimension $t$ over $\mathbf{F}_p$ with $ \mathbf{F}_p \subseteq V \subset k$ and 
$ \rho(u)(t)= t/(ut+1)$ for $ u \in V$ is the action of $V$ on $A\psl t \psr$. Also, if present, the factor ${\mathbf Z}/2= \langle w \rangle$ acts by multiplication by $-1$ on $t$. From \cite{CK}, 4.4.1, 4.4.2 and 4.4.5, and \cite{CM}, Proposition 3.2 and 3.7, after some simplification we find that the versal deformation ring is given by $$R= k\psl\alpha,x_2,\ldots,x_t\psr /(\alpha^{\frac{p-1}{2}},\alpha x_2,\ldots,\alpha x_t)$$ and the versal deformation is given by 
$$ \displaystyle \left\{ \begin{array}{lll} \xstrut \tilde{\rho}(1)(t)&=& \displaystyle \frac{t+\alpha}{t+\alpha+1},\\ \xstrut \tilde{\rho}(u_i)(t)&=& \displaystyle \frac{A_{u_i}t+B_{u_i}}{(C_{u_i}+x_i)t+D_{u_i}}, \quad 2 \leq i \leq t \\ \xstrut \tilde{\rho}(w)(t)&=& - t - \alpha,
\end{array} \right. $$ 
with the ``Formal Chebyshev Polynomials'' $$A_u=\sum_{j=0}^{\frac{p-1}{2}} \binom{u+j-1}{2j} \alpha^j, \quad C_u=\sum_{j=0}^{\frac{p-1}{2}-1} \binom{u+j}{2j+1} \alpha^j,$$ and $$ \quad B_u=\alpha C_u, \quad D_u=A_u+B_u. $$ 

We now prove universality. Assume that $A$ is an object of $\kate$ such that the map $\Phi_A \colon h_R(A) \to D_{\rho}(A)$ is not injective. Choose $\fhi_1, \fhi_2 \in h_{R}(A)$ with the same image in $D_{\rho}(A)$. Write $\alpha_1=\varphi_1(\alpha)$, $\alpha_2=\varphi_2(\alpha)$, $y_i=\varphi_1(x_i)$ and $z_i=\varphi_2(x_i)$. Since the corresponding functor for the action of the subgroup ${\mathbf F}_p \subseteq V$ has already been established to be universal in Proposition \ref{Z/p}, we get $\alpha_1=\alpha_2$. We write this element as $\alpha$, and now we only have to prove the other deformation parameters $y_i$ and $z_i$ coincide. 

By Propositions \ref{elab} and \ref{conjugation}, we can simultaneously conjugate the matrices in the image of $\fhi_1^* \circ \tilde{\rho}$ to the corresponding matrices in the image of $\fhi_2^* \circ \tilde{\rho}$ by a matrix $\gamma$ in $\PB_{A',k}$ for a quadratic extension $A \subseteq A'$. Now $\gamma$ commutes projectively with $M:=\fhi_1^* \circ \tilde{\rho}(1)=\fhi_2^* \circ \tilde{\rho}(1)= \big( \begin{smallmatrix} 1 & \alpha \\ 1 & \alpha+1 \end{smallmatrix} \big)$, a matrix with $\bar{M}$ not diagonal. Hence by Lemma \ref{commute}, $\gamma$ and $M$ are affinely dependent.
Since all $\tilde{\rho}(u_i)$ are affinely dependent on $M$ via $$\tilde{\rho}(u_i) = (C_u+x_i) M + (A_u-C_u-x_i) \cdot \mathrm{id}, $$
$\gamma$ also commutes with all (specializations of) $\tilde{\rho}(u_i)$, and therefore in the remaining conjugating equations $\gamma \fhi_1^* \circ \tilde{\rho}(u_i) = \mu_i \fhi_2^* \circ \tilde{\rho}(u_i) \gamma$ (for some invertible $\mu_i$) we can divide by $\gamma$ to get 
\begin{equation*} \label{zgliszcza} \mu_i  \Big(\begin{matrix} A&\alpha C \\ C+y_i&A+\alpha C
\end{matrix}\Big) =  \Big(\begin{matrix} A&\alpha C \\ C +z_i &A+\alpha C
\end{matrix}\Big)  \end{equation*} for some $A,C$  whose precise form is irrelevant for this proof, but we note $\bar{A}=1$. Hence the top left entry gives $\mu_i=1$, and from the lower left entry we indeed find that $y_i=z_i$. \qed
\end{se}

\begin{rem}\label{zadziornosc}
In \cite{CM}, the computation of the versal deformation ring in these cases depends heavily on Lemma 3.6 on page 245 of loc.\ cit.; the argument in the last sentence at the bottom of that page that allows one to conclude an equality of deformation parameters should be replaced by the universality of the versal deformation for ${\mathbf Z}/p$ from Proposition \ref{Z/p} in this paper. This correction does not create interdependencies of proofs.  
\end{rem}

\begin{se}[$G=({\mathbf Z}/2)^t, p=2, t \geq 3$] Let $G=(\mathbf{Z}/2)^t$ with $t \geq 3$. Write $G=V$ with $V$ a sub-$\mathbf{F}_2$-vector space of $k$ with basis $\{1,u_2,\ldots,u_t\}$. We can suppose that $\rho \colon G \to \Gamma_k$ is given by $\rho(u)(t)=t/(ut+1)$ for $u\in V$. The computation of the versal deformation on page 454 of \cite{CK} is false, and should be replaced by: the versal deformation ring of $\rho$ is given by $$R=k\psl\alpha,x_3,\ldots,x_t\psr$$ and the versal deformation is given \emph{on generators} by $$ \left\{ \begin{array}{llll} \xstrut \tilde{\rho}(1)(t)&=&  \displaystyle \frac{t+\alpha}{t+1},\\ \xstrut \tilde{\rho}(u_2)(t)&=&  \displaystyle \frac{t+\alpha u_2}{u_2t+1}, \\  \xstrut \tilde{\rho}(u_i)(t)&=&\displaystyle  \frac{t+\alpha (u_i+x_i)}{(u_i+x_i)t+1}, \quad 3 \leq i \leq t.\end{array} \right. $$ As for the proof of versality, these matrices are easily seen to satisfy the relations of the generators of $G$ without any conditions on the deformation parameters, and since the resulting deformation ring is smooth and the induced map on tangent spaces is an isomorphism by the calculation of group cohomology  in loc.\ cit., we are done.

We now prove universality. The proof is very similar to the previous case, though the actual versal deformations are different. Assume that $A$ is an object of $\kate$ such that the map $\Phi_A \colon h_R(A) \to D_{\rho}(A)$ is not injective. Choose $\fhi_1, \fhi_2 \in h_{R}(A)$ with the same image in $D_{\rho}(A)$. Write $\alpha_1=\varphi_1(\alpha)$, $\alpha_2=\varphi_2(\alpha)$, $y_i=\varphi_1(x_i)$ and $z_i=\varphi_2(x_i)$. Since the corresponding functor for the equicharacterstic deformation of the action of the subgroup ${\mathbf F}_2(u_2) \subseteq V$ has already been established to be universal in Proposition \ref{smut}, we get $\alpha_1=\alpha_2$. We write this element as $\alpha$, and now we only have to prove the other deformation parameters $y_i$ and $z_i$ coincide. 

Observe that the image $\tilde{\rho}(G)$ is a set of affinely dependent matrices via $$\tilde{\rho}(u_i)=(u_i+x_i)\tilde{\rho}(1)+(u_i+x_i+1) \cdot \mathrm{id}.$$ By Propositions \ref{kar2} and \ref{conjugation}, we can simultaneously conjugate the matrices in the image of $\fhi_1^* \circ \tilde{\rho}$ to the corresponding matrices in the image of $\fhi_2^* \circ \tilde{\rho}$ by a matrix $\gamma$ in $\PB(A')$ for a quadratic extension $A \subseteq A'$. Now $\gamma$ commutes projectively with both $\fhi_1^* \circ \tilde{\rho}(1)=\fhi_2^* \circ \tilde{\rho}(1)= M:= \big( \begin{smallmatrix} 1 & \alpha \\ 1 & 1 \end{smallmatrix} \big)$, and $\fhi_1^* \circ \tilde{\rho}(u)=\fhi_2^* \circ \tilde{\rho}(u)$ and by Remark \ref{qual}, this implies an affine dependency between $M$ and $\gamma$. Since all $\tilde{\rho}(u_i)$ are affinely dependent on $M$, 
$\gamma$ also commutes with all (specializations of) $\tilde{\rho}(u_i)$, and therefore in the remaining conjugating equations $\gamma \fhi_1^* \circ \tilde{\rho}(u_i) = \mu_i \fhi_2^* \circ \tilde{\rho}(u_i) \gamma$ (for some invertible $\mu_i$) we can divide by $\gamma$ to get 
\begin{equation*} \mu_i \Big(\begin{matrix} 1&\alpha (u_i+y_i) \\ u_i+y_i&1
\end{matrix}\Big) =  \Big(\begin{matrix} 1&\alpha (u_i+z_i) \\ u_i +z_i &1
\end{matrix}\Big) \end{equation*} for some  invertible $\mu_i \in (A')^*$. Hence we indeed find $\mu_i=1$ from the top left entry, and then $y_i=z_i$ from the lower left entry. \qed
\end{se}

\begin{remark}
The proof of Proposition 3.8(ii) on page 250 of \cite{CM} should be replaced by a direct calculation that, however, leads to the same result. More precisely, one assumes there is a lift to a ring $A$ of characteristic $4$ and considers its reduction to $A/2\maks_A$. Write the lifts of three generators explicitly as the universal equicharacteristic lifts from \cite{CK} plus $2$ times an unknown power series in $k\psl t \psr$. The relations between the generators (having order two and commuting) give rise to a system of linear equations in $k$ between the coefficients of the first three terms of those three power series, that has no solutions. 
\end{remark}

\begin{se}[$G=({\mathbf Z}/p)^t \rtimes {\mathbf Z}/m \mbox{ with } m \geq 3 \mbox{ and either } t \geq 2 \mbox{ or } p \geq 3$] Set $P=({\mathbf Z}/p)^t$ and $G=P \rtimes {\mathbf Z}/m$; the previous results imply that (with obvious notations) $D_P$ is pro-representable. There is a  restriction map $D_G \rightarrow D_P$, and we have a surjective map $R_P \rightarrow R_G$ from looking at the explicit form of those rings in \cite{CM} (or by a general theory developed in \cite{Kuba}). Therefore, in the commutative diagram 
$$\xymatrix{h_{R_G}\ar[r]\ar[d] & D_G\ar[d]\\ h_{R_P}\ar[r] & D_P}$$ the left vertical map is injective and the bottom map an isomorphism, hence the top map is also injective. \qed
\end{se}

Since we have now treated all possible weakly ramified actions of a finite group in $\Gamma_k$, this finishes the proof of the two main theorems from the introduction. \qed

\section{Final remarks}

\begin{se}
We remark on the general strategy employed here to show pro-representability, and compare it to other instances of proofs of pro-representability. 

We start off with a versal deformation over a versal ring $R$ and assume pro-re\-pre\-sen\-tabi\-li\-ty fails at some ring $A \in \kate$: there are two different morphisms from $R$ to $A$ that induce the same deformation, that is, the same representation up to $\Gamma_A$-conjugacy. Our study of the Nottingham group over $\kate$ shows that we can conjugate the image of one of these representations into $\Gamma^0_{A'}$, but only for an extension $A \subseteq A'$ that depends on the versal deformation: we have to extract a square root. Since the image is now in $\Gamma^0_{A'}$, we can change the conjugation in $\Gamma_{A,k}$ into a ``matrix conjugation'' in $\PB_{A',k}$. This, however, leads to some trace and determinant equalities that show the original morphisms to be equal, contradiction. 

It is important to note that the proof doesn't show that any $\Gamma_{A,k}$-conjugacy can be replaced by a $\PB_{A,k}$-conjugacy. 
\end{se}

\begin{se}
The classical proofs of pro-representability in case of the deformation theory of a curve of genus $g \geq 2$ (\cite{Sernesi}, \cite{Schlessinger}), or a group action on one of those (\cite{BM}), or of an absolutely irreducible linear representation (\cite{Mazur}) can be abstractly viewed as follows: 

\begin{prop} Let $F_0 \colon \kate \rightarrow \mathsf{Sets} $ denote a pro-representable, not necessarily finite dimensional  functor, and let ${\mathfrak G} \colon \kate \rightarrow \mathsf{Groups}$ denote a pro-representable group functor that acts on $F_0$. Then the functor $F:=\mathfrak{G} \backslash F_0$ is pro-representable if the following two properties hold:
\begin{enumerate}
\item[{\bf P1}] for any surjection $A \twoheadrightarrow A_0$ in $\kate$, ${\mathfrak G}(A) \twoheadrightarrow {\mathfrak G}(A_0)$ is surjective;

\item[{\bf P2}] for any surjection $A \twoheadrightarrow A_0$ in $\kate$, any $\xi \in F_0(A)$ and $g_0 \in {\mathfrak G}(A_0)$ such that $$g_0 \xi {\big|}_{A_0} = \xi{\big|}_{A_0},$$ we can find $g \in {\mathfrak G}(A)$ with $g{\big|}_{A_0}=g_0$ and $g \xi = \xi$. \qed
\end{enumerate}
\end{prop} 

Since $F_0$ can have an infinite-dimensional tangent space, the proof uses Grothen\-dieck's original criterion for pro-representability (\cite{Groth}) but is easy and further left to the reader. 

In the above applications, $F_0$ is the functor taking $A \in \kate$ to the set of all lifts to $A$ of the object under consideration  (curve, group action, representation), and $\mathfrak G$ the functor describing the equivalence relation that is divided out in constructing the actual deformation functor (infinitesimal automorphisms, linear conjugacy, \dots). For example, for an absolutely irreducible group representation, by Schur's Lemma the commutator of the image of the linear representation is the set of scalar matrices, and those are seen to lift in the sense of {\bf P2} (\cite{Mazur}).  In our setting of local representations, {\bf P1} is clear, so if one could establish directly that {\bf P2} hold for $F_0(A)$ the set of lift of a given representation and ${\mathfrak G}(A)=\Gamma_{A,k}$ acting by conjugation, then pro-representability would follow. Also, one may wonder whether it is exactly for the non-pro-representable cases in our main theorems that {\bf P2} fails. These appear to be difficult problems concerning power series over Artin rings. Here, we only show that one can reduce in the weakly ramified case to the study of conjugation by $\PB_{A',k}$ over varying extensions $A'$.  
\end{se}

\begin{se}
In \cite{Dokchitserarxiv}, p.4, Tim Dokchitser has asked whether if $D$ is a not necessarily pro-representable functor but admits a hull $h_R$, there exists a  group functor $\mathfrak G$ acting on $h_R$ such that $D \cong \mathfrak{G} \backslash h_R$. The easiest test case that comes out of our analysis is the deformation functor $D$ of a weakly ramified local representation $\rho \colon {\mathbf Z}/2 \to \Gamma_k$ over a field $k$ of characteristic $2$: $D$ has a versal hull $R=W(k)\psl \alpha \psr$, but is not pro-representable. 

First note the following easy observation: 
\begin{prop} $D$ is \emph{not} a quotient of a pro-representable functor by a constant group action.
\end{prop}
\begin{proof}  If $D$ is a quotient of $h_R$ by a constant group functor, it takes injections to injections, hence it suffices to remark that this is not the case for $D$ (note: in \cite{Dokchitserarxiv}, Theorem 1, one even finds \emph{necessary and sufficient} conditions for $D$ to be a quotient by a constant group functor).
Take $A=k[\varepsilon]/\varepsilon^2$ and set $$m= \Big( \begin{matrix} 1 & 0 \\ 1 & 1 \end{matrix} \Big) \mbox{ and } n= \Big( \begin{matrix} 1 & \varepsilon \\ 1 & 1 \end{matrix} \Big).$$ Take $$B=k[\varepsilon,w]/[w^2-\varepsilon,w^3] \simeq k[w]/w^3.$$ From Proposition \ref{psychopatologia}, it follows that $m$ and $n$ define different elements of $D(A)$, but define the same element of $D(B)$; since in $\PGL(2,B)$, $n$ has a representative of determinant 1  (namely, $(1+w)^{-1} n$). Hence the map $D(A) \to D(B)$ is not injective.
\end{proof}

Though we cannot answer the question of Dokchitser for this group action, our techniques can be used to squeeze the functor $D$ between two group functor quotients of $h_R$, as follows. The versal deformation is given in \ref{two}, so let $\alpha, \beta \in A$ in $\kate$ with $m_\alpha$ and $m_\beta$ conjugate in $\Gamma_{A,k}$, where we put $m_\alpha:= \big( \begin{smallmatrix} 1 & \alpha \\ 1 & -1 \end{smallmatrix} \big)$. Our results imply that for $A'=A[{\mathfrak z}]$ and ${\mathfrak z}^2-2{\mathfrak z} -\alpha=0$, there exists a matrix $\gamma \in \PB_{A',k}$ and a constant $\lambda \in (A')^*$ with $\bar{\lambda}=1$  and $\gamma m_\alpha = \lambda m_\beta \gamma$. 

We now show how to make this equivalence under $\PB_{A',k}$ into a group functor action on $h_R$. This requires some work since $\PB_{A',k}$ itself has $A'$ varying with $A$. However, we know by Proposition \ref{psychopatologia} that the equivalence occurs exactly when $\lambda^2 = \frac{1+\alpha}{1+\beta}$. Writing $\lambda=\lambda_1 + \lambda_2 {\mathfrak z}$, we get $$\lambda^2 = \lambda_1^2 + \alpha \lambda_2^2  +2 \lambda_2 (\lambda_1+\lambda_2){\mathfrak z} \in A,$$ and since $\{1,{\mathfrak z}\}$ are independent over $A$, we find  $\lambda^2 = \lambda_1^2 + \alpha \lambda_2^2$ with $2 \lambda_2(\lambda_1+\lambda_2)=0$.   We switch variables to $(a,b):=(\lambda_1 + \lambda_2,\lambda_2)$, and the condition becomes $2ab=0$, so $\lambda=(a-b)^2+\alpha b^2 = a^2 + b^2(\alpha+1)$. We also still have the condition $\bar{\lambda_1}=1$ that translates to $\bar{a}+\bar{b}=1$.  We want to see whether we can make the set  $${\mathfrak M}(A):= \{ (a,b) \in A \times A \colon \bar{a}+\bar{b}=1 \, \wedge 2ab=0 \}$$ act on the transformations $m_\alpha$ (or equivalently, on $\alpha$) such that equivalent matrices are in the same orbit. We therefore investigate whether we can impose on the set ${\mathfrak M}(A)$ a binary operation $\times$ with an action on $\alpha$ that describes the given operation 
\begin{equation}\label{ast} \alpha \mapsto \beta= (a,b) \ast \alpha := \frac{1+\alpha}{a^2+b^2 (\alpha+1)}-1. \end{equation} We find by direct calculation that 
$$ (a,b) \ast ((c,d) \ast \alpha) = \frac{1 + \alpha}{a^2c^2+(a^2d^2+b^2)(\alpha+1)}-1. $$
Since  $2ab=0$, we have $a^2d^2+b^2=(ad+b)^2$, hence 
we should put
\begin{equation} \label{times} (a,b) \times (c,d) := (ac, ad+b). \end{equation}
This makes $({\mathfrak M}, \times)$ into  a monoid functor ${\mathfrak M} \colon \kate \rightarrow {\mathsf{Monoids}}$. 
One indeed checks this operation is associative with neutral element $(1,0)$. The inverse of $(a,b)$ is $(a^{-1},-ba^{-1})$, which only exists if $a \in A^*$. Hence this monoid functor has a group subfunctor ${\mathfrak G} \colon \kate \rightarrow {\mathsf{Groups}}$ given by \begin{equation}\label{G} {\mathfrak G}(A):=  \maks[2](A) \rtimes U(A), \end{equation} with $\maks[2](A) = \{ b \in \maks_A \colon 2b=0 \}$ and $U(A)=1+\maks_A$ and the operation $\times$ on $\mathfrak G$ is induced from the usual Borel semidirect product structure on ${\mathbf G}_a \rtimes {\mathbf G}_m$ of which ${\mathfrak G}(A)$ is a subgroup. 
 
The monoid $\mathfrak M$ and its subgroup $\mathfrak G$ act on $h_R$, since $m \ast (n \ast \alpha) = (m \times n) \ast \alpha$ for any $m,n \in {\mathfrak M}(A)$ by construction. Furthermore, $\mathfrak M$ acts on $h_R$ in such a way that we get a morphism of functors $D \rightarrow {\mathfrak M} \backslash h_R$. We now show that $D$ even maps to the group functor quotient $D \rightarrow {\mathfrak G} \backslash h_R.$ Indeed, suppose two elements $\alpha$ and $\beta$ are connected by $\alpha=m\beta$ for $m$ and $\beta=n\alpha$ with $m=(m_1,m_2)$ and $n=(n_1,n_2)$  in ${\mathfrak M}(A)$ but not in ${\mathfrak G}(A)$. Then $\alpha$ is fixed under the action of $mn$. With the explicit formula above, this gives $$m^2_1 n^2_1 + (m_1 n_2 + m_2)^2 ( 1 +\alpha)=1$$ and hence 
$1+\alpha$ is a square, namely, the square of $(1+m_1 n_1)/(m_1 n_2+m_2)$, since $2m_1n_1=0$ (note that $m_1 n_2 + m_2$ is invertible since we have chosen $m$ outside ${\mathfrak G}(A)$). But then it follows immediately that $m_\alpha$ is $\PB_{A,k}$-conjugate to $m_\beta$, hence there exists $g \in {\mathfrak G}$ with $\alpha=g\beta$, as was to be shown. 

Thus, we get 
\begin{prop}Let $R=W(k)\psl \alpha \psr$ denote the versal deformation ring of a weakly ramified involution $t \mapsto 1/(1+t)$ over a field $k$ of characteristic $2$ with versal deformation $t \mapsto (t+\alpha)/(t-1)$. Then the functor $h_R$ admits an action of group functor ${\mathfrak G} = \maks[2] \rtimes U$ given in \textup{(\ref{ast})} via the representation \textup{(\ref{G})} and \textup{(\ref{times})} such that 
there exists a morphism of functors  $D \rightarrow {\mathfrak G} \backslash h_R.$ \qed
\end{prop}
If, on the other hand, we consider the group subfunctor $U = \{(-,0) \in {\mathfrak M}\}$ acting on $h_R$ by conjugation, we find that the $\PB_{A',k}$-conjugacy that results to be defined already over $A$, 
 and thus we get a morphism of functors in the other direction $U \backslash h_R \rightarrow D.$ 
We arrive at 
\begin{prop} There exists a sequence of morphisms of functors
$$h_R \rightarrow U \backslash h_R \xrightarrow{\Psi} D \xrightarrow{\Phi} (\maks[2] \rtimes U) \backslash h_R. \qed $$ 
\end{prop}
At present, we do not know whether any of these morphisms is an isomorphism. For $\Psi$ to be an isomorphism means that any $\Gamma_{A,k}$-conjugacy can be transformed into a $\PB_{A,k}$-conjugacy. For $\Phi$ to be an isomorphism means the following: set $A_\alpha:=A[{\mathfrak z}]$ with ${\mathfrak z}^2=2{\mathfrak z}+\alpha$; then if $m_\alpha$ and $m_\beta$ are conjugate by an element of  $\PB_{A_\alpha,k}$ and by an element of $\PB_{A_\beta,k}$, they are $\Gamma_{A,k}$-conjugate. Again, these appear to be difficult questions concerning power series over Artin rings. It is precisely trying to circumvent these difficulties in the general case that led to our method. Note that any of these two maps being an isomorphism would confirm Dokchitser's conjecture for a weak local involution in characteristic two.

\end{se}

\begin{small}
\bibliographystyle{amsplain}
\providecommand{\bysame}{\leavevmode\hbox to3em{\hrulefill}\thinspace}
\providecommand{\MR}{\relax\ifhmode\unskip\space\fi MR }
\providecommand{\MRhref}[2]{%
  \href{http://www.ams.org/mathscinet-getitem?mr=#1}{#2}
}
\providecommand{\href}[2]{#2}

\end{small}

\end{document}